\documentclass[11pt]{article}
\usepackage{epsfig,amsmath,amsfonts,latexsym,subfigure,amssymb}
\usepackage[latin1]{inputenc}
\usepackage{graphicx,epic}
\usepackage[all]{xy}
\newtheorem{theorem}{Theorem}[section]

\newtheorem{lemma}[theorem]{Lemma}
\newtheorem{prop}[theorem]{Proposition}
\newtheorem{corr}[theorem]{Corollary}

\def\0{{\bf 0}}
\def\P{{\bf P}}
\def\E{{\bf E}}

\begin{document}

\title{Epidemics on random graphs with tunable clustering}

\author{Tom Britton \thanks{Department of Mathematics, Stockholm
University, 106 91 Stockholm, Sweden. E-mail: {\tt tomb@math.su.se, mia@math.su.se, andreas@math.su.se, lindholm@math.su.se}} \and Maria Deijfen$^*$ \and Andreas N.\ Lager\aa s$^*$ \and Mathias Lindholm$^*$}

\date{July 2007}

\maketitle

\thispagestyle{empty}

\begin{abstract}

\noindent In this paper, a branching process approximation for the
spread of a Reed-Frost epidemic on a network with tunable
clustering is derived. The approximation gives rise to expressions
for the epidemic threshold and the probability of a large outbreak in the epidemic. It is investigated how these quantities varies
with the clustering in the graph and it turns out for instance
that, as the clustering increases, the epidemic threshold
decreases. The network is modelled by a random intersection graph,
in which individuals are independently members of a number of
groups and two individuals are linked to each other if and only if
they share at least one group.

\vspace{0.5cm}


\vspace{0.5cm}

\end{abstract}

\section{Introduction}

This paper is concerned with Reed-Frost epidemics modified to take
place on random networks. Introduced in 1928 by two medical
researchers, Lowell Reed and Wade Frost, the Reed-Frost model is
one of the simplest stochastic epidemic models. The spread of the
infection takes place in generations: Each individual that is
infective at time $t$ ($t=0,1,\ldots$) independently makes
contacts with all other individuals in the population with some
probability $p$, and if a contacted individual is susceptible, it
becomes infected at time $t+1$. Also at time $t+1$, the infective
individuals from time $t$ are removed from the epidemic process.

The behavior of the Reed-Frost model is well understood, see e.g.\
von Bahr and Martin-L\"{o}f (1980). A crucial assumption which
simplifies the analysis of the model is that the population in
which the epidemic takes place is taken to be homogeneously
mixing, that is, an infective individual is assumed to make
contacts with all other individuals in the population with the
same probability. This assumption is of course very unrealistic,
since, in a real-life epidemic, an infective individual is much
more likely to infect individuals with whom he/she has some kind
of social connection. The Reed-Frost model can easily be adapted
to incorporate this type of heterogeneity by introducing a graph
to represent the social structure in the population and then
stipulating that infective individuals can only infect their
neighbors in the social network; see Section 3. This modification
makes the analysis of the model two-fold. Firstly, one wants to
find a realistic model for the underlying social network, and,
secondly, one wants to study the behavior of the epidemic on this
graph.

Large complex networks such as social contact structures, the
internet and various types of collaboration networks have received
a lot of attention during the last few years; see e.g.\
Dorogovtsev and Mendes (2003) and Newman et al.\ (2006) and the
references therein. As for social networks, one of their most
striking features is that they are highly clustered, meaning
roughly that there is a large number of triangles and other short
cycles; see e.g.\ Newman (2003). This is a consequence of the fact
that friendship circles are typically strongly overlapping so that
many of our friends are also friends of each other. A model that
captures this in a natural way is the so called random
intersection graph, which is described in Section 2. Roughly, the
idea of the model is that people are members of groups ---
families, schools, workplaces etc.\ --- and an edge is drawn
between two individuals if they share at least one group. If the
relation between the number of individuals and the number of
groups is chosen appropriately, this leads to a graph where the
amount of clustering can be tuned by adjusting the parameters of
the model.

An important goal of network modelling is to investigate how the
structure of the network affects the behavior of various types of
dynamic processes on the network; see Durrett (2006) for an
overview. When it comes to epidemics, Andersson (1999) is a comprehensible
introduction, in which expressions for the epidemic threshold, the
probability of a large outbreak and the final size of the epidemic
are derived in a heuristic way for a number of underlying graphs.
Here, the epidemic threshold, commonly denoted by $R_0$, is
defined as a function of the parameters of the model such that a
large outbreak in the epidemic has positive probability if and
only if $R_0>1$. In epidemic modelling, a common technique for
deriving expressions for the epidemic threshold and the
probability of a large outbreak is to use branching process
approximations of the early stages of the epidemic. However, when
studying epidemics on networks, dependencies between the edges in
the graph tend to make branching process approximations more
complicated. Results for epidemics on graphs with arbitrary degree
distribution can be found in Andersson (1998), and
Erd\H{o}s-R\'{e}nyi graphs and some extensions thereof are dealt
with in Neal (2004, 2006). There is however very little work done
on more complicated graph structures.

The aim of this paper is to give a rigorous analysis of how
clustering in a network affects the spread of an epidemic. The
network is modelled by a random intersection graph with tunable
clustering and we then let a Reed-Frost epidemic propagate on this
graph. Comparing the epidemic with a certain branching process
yields (implicit) expressions for the epidemic threshold and the
probability of a large outbreak. Numerical evaluations reveal
that, as the clustering increases, the epidemic threshold
decreases --- that is, large outbreaks are possible for larger
parts of the parameter space --- but also that the actual value of
the probability of a large outbreak decreases as the clustering
approaches its maximal value. To our knowledge, this is the first
rigorous investigation of how the spread of an epidemic is
affected by clustering.

In Newman (2003:2), the effect of clustering on 
epidemics is studied by heuristic means, and calculations therein
indicate that indeed the epidemic threshold should decrease as the
clustering increase. Furthermore, Trapman (2007) studies epidemics
on graphs with a given expected number of triangles, but the
construction of the graph is more involved there. We mention also
the work by Ball et al.\ (1997) on the so called household model,
which describes the spread of an epidemic in a population with
group structure. The model there however is not formulated in
terms of an underlying graph and the concept of clustering is not
considered.

The paper is organized as follows. In Section 2, random
intersection graphs and their properties are described in more
detail. Section 3 contains the main result --- a comparison of a
Reed-Frost epidemic on a random intersection graph with a
branching process --- and its proof. In Section 4, the final size
of the epidemic is commented on. It is observed that a thinned
random intersection graph is in fact not a random intersection
graph, implying that results concerning the component structure in
a random intersection graph cannot be used to draw conclusions
about the final size of the epidemic. In Section 5, numerical
results are presented and, finally, Section 6 contains a short
discussion.

\section{Random intersection graphs}

Random intersection graphs were introduced in Singer (1995) and
Karo\'{n}ski et al.\ (1999). In its simplest form, the model is
defined as follows: Given a set $\mathcal{V}$ of $n$ vertices and a
set $\mathcal{A}$ of $m$ auxiliary vertices, construct a bipartite
graph $\mathcal{B}_{n,m,r}$ by letting each edge between vertices
$v\in\mathcal{V}$ and $a\in\mathcal{A}$ exist independently with
probability $r$. The random intersection graph $\mathcal{G}_{n,m,r}$
with vertex set $\mathcal{V}$ is obtained by connecting two vertices
$v,w\in\mathcal{V}$ if and only if there is a vertex
$a\in\mathcal{A}$ such that $a$ is linked to both $v$ and $w$ in
$\mathcal{B}_{n,m,r}$. This construction can be generalized in
various ways --- see e.g.\ Godehardt and Jaworski (2002) and Deijfen
and Kets (2007) --- but in this paper we will stick to the above
formulation. We will also specialize to the case when
$m=\lfloor\beta n^\alpha \rfloor$ for some constants
$\alpha,\beta>0$; see Karo\'{n}ski et\ al.\ (1999) for a motivation
of this choice of $m$. In fact, to get a graph with tunable
clustering, we will soon take $\alpha=1$.

If the vertices in $\mathcal{V}$ and $\mathcal{A}$ are thought of as
individuals and groups respectively, then the random intersection
graph provides a model for a social network where individuals are
connected if there is at least one group where they are both
members. The probability that two individuals do not share any group
is $(1-r^2)^m$, implying that the edge probability in the random
intersection graph is $1-(1-r^2)^m$, and hence the expected degree
of a fixed vertex is
$$
(n-1)(1-(1-r^2)^m)=\beta r^2n^{1+\alpha}+O(r^4n^{1+2\alpha}).
$$
To keep this expression bounded as $n\to\infty$, we let $r=\gamma
n^{-(1+\alpha)/2}$ for some $\gamma>0$. The expected degree then
tends to $\beta\gamma^2$ as $n\to\infty$.

As for the asymptotic distribution of the vertex degree with the
above choices of $m$ and $r$, it is shown in Stark (2004) to be a
point mass at 0 for $\alpha<1$, a compound Poisson distribution,
describing the law of a sum of a Poisson($\beta\gamma$)
distributed number of independent Poisson($\gamma$) variables for
$\alpha=1$, and a Poisson($\beta\gamma^2$) distribution for
$\alpha>1$. To see this, note that the number of groups that an
individual belongs to is binomially distributed with mean
$mr=\beta\gamma n^{(\alpha-1)/2}$. For $\alpha<1$, this goes to 0
as $n\to\infty$, explaining the point mass at 0. For $\alpha=1$,
the number of group memberships per individual is asymptotically
Poisson($\beta\gamma$) distributed, and the sizes of the groups
are Poisson($\gamma$), with overlaps between groups being very
unlikely if $n$ is large, indicating that the degree distribution
should indeed be compound Poisson. When $\alpha>1$, each
individual belongs to infinitely many groups as $n\to\infty$. This
means that the edge indicators are asymptotically independent,
which suggests a Poisson distribution for the vertex degree. In
fact, for $\alpha>1$, the random intersection graph is similar to
the standard Erd\H{o}s-R\'{e}nyi random graph; see Fill et al.\
(2000).

Moving on to the clustering in the graph, for two vertices
$v,w\in\mathcal{V}$, let $I_{vw}$ denote the edge indicator for the
edge between $v$ and $w$ in $\mathcal{G}_{n,m,r}$, and write $\P_n$
for the probability measure of $\mathcal{G}_{n,m,r}$. We then define
the clustering as
$$
c=c_{\alpha,\beta,\gamma}:=\lim_{n\to\infty}\P_n(I_{vw}=1|I_{vu}I_{wu}=1),
$$
that is, $c$ is the limiting conditional probability that there is
an edge between two vertices $v$ and $w$, given that they have a
common neighbor $u$. The expected number of groups that individual
$u$ belongs to is $\beta\gamma n^{(\alpha-1)/2}$, which goes to 0,
$\beta\gamma$ or infinity, depending on whether $\alpha<1$,
$\alpha=1$ or $\alpha>1$. As a consequence, the limiting
probability that two individuals $v$ and $w$ who both share a
group with $u$ in fact share the \emph{same} group with $u$ --- thus being
connected to each other --- will behave differently depending on
the value of $\alpha$. More specifically, it is shown in Deijfen
and Kets (2007) that
$$
c_{\alpha,\beta,\gamma} = \left\{
\begin{array}{ll}
1 & \mbox{if } \alpha<1;\\
(1+\beta\gamma)^{-1} & \mbox{if } \alpha=1;\\
0 & \mbox{if } \alpha>1.
\end{array}\right.
$$
In view of the result from Stark (2004) concerning the degree
distribution and the characterization of the clustering from
Deijfen and Kets (2007), the best choice if we want to use a random intersection graph to describe a social network seems to be
$\alpha=1$. This gives rise to a model where both the mean degree
and the clustering can be tuned by adjusting the parameters
$\beta$ and $\gamma$. More precisely, with $D$ denoting the
limiting degree of a fixed vertex, we have that
$$
\E[D]=\beta\gamma^2\quad\mbox{and}\quad c=(1+\beta\gamma)^{-1}.
$$
For the remainder of this paper we fix $\alpha=1$ and write
$\mathcal{G}_{\beta,\gamma}^{(n)}=\mathcal{G}^{(n)}$ and
$\mathcal{B}_{\beta,\gamma}^{(n)}=\mathcal{B}^{(n)}$ for the
corresponding random intersection graph and its underlying
bipartite graph (omitting the subscripts when the dependence on
$\beta$ and $\gamma$ does not need to be emphasized).

\section{The epidemic model and an approximating branching process}

Consider a closed homogeneous population consisting of $n$
individuals, labelled $v_1,\ldots,v_n$, with a social structure
represented by a random intersection graph $\mathcal{G}^{(n)}$. We
will use the Reed-Frost dynamics to describe the spread of an
infection in this population. The social graph $\mathcal{G}^{(n)}$
is assumed to be fixed throughout the spread of the infection.
Furthermore, for simplicity, we start with one single randomly
selected infective individual at time 0, the rest of the
population being susceptible. Without loss of generality, we
assume that the initial infective, which will be referred to as
the index case, is individual $v_1$. An individual that is
infective at time $t$ ($t=0,1,\ldots$) contacts each one of its
neighbors in $\mathcal{G}^{(n)}$ independently with some
probability $p$, and if a contacted neighbor is susceptible, it
becomes infective at time $t+1$. The individuals that were
infective at time $t$ are removed from the epidemic process at
time $t+1$ (by immunity or death) and take no further part in the
spread of the infection.

We will be concerned with the set $\mathcal{E}^{(n)}$ of
individuals that are ultimately affected by the above epidemic.
More precisely, we will construct a branching process that can be
used to determine whether $\mathcal{E}^{(n)}$ is finite or
infinite in the limit as $n\to\infty$. To this end, first note
that $\mathcal{E}^{(n)}$ can be identified with the cluster
containing the index case in an edge percolation process on
$\mathcal{G}^{(n)}$ in which each edge is open independently with
probability $p$. Open edges in the percolation process are
interpreted as possible transmission links for the disease, that
is, if one of the vertices of an edge is infective at time $t$ and
the other one is not, then the uninfected vertex becomes infective
at time $t+1$. Furthermore, if we consider the percolation cluster
of a particular vertex restricted to a subgraph of
$\mathcal{G}^{(n)}$, then the size of this cluster has the same
distribution as the final size of a Reed-Frost epidemic on the
subgraph. In particular, if $R_k$ is the size of the percolation
cluster of a given vertex belonging to a complete subgraph with
$k$ nodes, then the distribution of $R_k$, denoted by $F_k$, is
that of a Reed-Frost epidemic initiated by one single individual
in a homogeneously mixing population of size $k$. The distribution
$F_k$ can be computed recursively; see Andersson and Britton
(2000: Section 1.2).

We now define the branching process that will be used to
approximate the epidemic process. To begin with, note that the
groups in a random intersection graph generate complete subgraphs
with sizes that are asymptotically Poisson($\gamma$) distributed.
Hence, in the limit as $n\to\infty$, the size of an outbreak
started by a given individual in a given group is

\begin{equation}\label{eq:R}
R\sim F:=\sum_{k=0}^{\infty}F_k\frac{\gamma^k}{k!}e^{-\gamma}.
\end{equation}

\noindent Recall that the number of groups that a given individual
is a member of is asymptotically Poisson($\beta\gamma$)
distributed. Let $f$ be the generating function of a sum of a
Poisson($\beta\gamma$) number of i.i.d.\ variables
$R_1,R_2,\dots$, all distributed as $R$, that is,

\begin{equation}\label{eq:f}
f(s) = \exp\left\{\beta\gamma(\E[s^R]-1)\right\},
\end{equation}

\noindent and let $\{Z(t):t\geq 0\}$ be a discrete time branching
process with offspring generating function $f$, that is,
$\E[s^{Z(1)}]=f(s)$. Finally, write $E$ for the total progeny in
such a process, that is,
$$
E=\sum_{t=0}^\infty Z(t).
$$

Let $E^{(n)}=|\mathcal{E}^{(n)}|$ denote the final size of a
Reed-Frost epidemic on a random intersection graph
$\mathcal{G}^{(n)}$. Our main result is the following theorem,
which will be proved by relating the initial phases of the
epidemic to a branching process with the same distribution as
$\{Z(t):t\geq 0\}$ as $n\to\infty$.

\begin{theorem}\label{main}
As $n\to\infty$, we have that $E^{(n)}\to E$ in distribution.
\end{theorem}

Define $\rho$ to be the smallest non-negative root of the equation
$f(\rho)=\rho$. It follows from standard results in branching
processes theory that $P(E=\infty)=1-\rho$ and that $\rho<1$ if
and only if $\E[Z(1)]>1$; see e.g.\ Athreya and Ney (1972).
Combining this with Theorem \ref{main} gives the following
corollary concerning the asymptotic behavior of the epidemic.

\begin{corr}\label{large_outbreak} Define $R_0:=\E[Z(1)]=\beta \gamma\E[R]$ and write $\pi=1-\rho$.
As $n\to\infty$, we have that
\begin{itemize}
\item[\rm{(a)}] $E^{(n)}\to\infty$ with probability $\pi$;
\item[\rm{\rm{(b)}}] $\pi>0$ if and only if $R_0>1$.
\end{itemize}
\end{corr}

Before we continue with the proof of the theorem, we will state
and prove a lemma concerning the bipartite graph
$\mathcal{B}^{(n)}$. To this end, for an arbitrary graph
$\mathcal{G}$ with vertex set $\mathcal{W}$, the subgraph of
$\mathcal{G}$ induced by some subset $\mathcal{W}'\subset
\mathcal{W}$ is defined to be the subgraph consisting the vertices
in $\mathcal{W}'$ together with all edges in $\mathcal{G}$ that
run between vertices in $\mathcal{W}'$. Let $\mathcal{C}^{(n)}(t)$
be the vertices of $\mathcal{B}^{(n)}$ at distance $t$ from vertex
$v_1$. Note that a vertex in $\mathcal{B}^{(n)}$ may be either an
individual (that is, a vertex $v\in\mathcal{V}$) or a group (that
is, an auxiliary vertex $a\in\mathcal{A}$ ), and that vertices at
odd distance from $v_1$ correspond to groups and vertices at even
distance to individuals.

\begin{lemma}\label{fodelsedag}
Let $\kappa>0$ be such that $1/\kappa>2\log(\beta\gamma^2)$. As
$n\to\infty$, the probability that the subgraph of
$\mathcal{B}^{(n)}$ induced by  $\mathcal{C}^{(n)}(\lfloor
\kappa\log n\rfloor)$, is a tree, tends to 1.
\end{lemma}

\noindent\textbf{Proof of Lemma \ref{fodelsedag}} We will build up
$\mathcal{C}^{(n)}(t)$ by a sequence $\{\mathcal{D}^{(n)}(t):t\geq
0\}$, constructed in such a way that
$\mathcal{C}^{(n)}(t)=\cup_{0\leq s \leq t}\mathcal{D}^{(n)}(s)$.
For odd $t$, the set $\mathcal{D}^{(n)}(t)$ will consist of groups
and for even $t$ by individuals. To begin with, by definition, we
have $\mathcal{C}^{(n)}(0)=\{v_1\}$, so necessarily
$\mathcal{D}^{(n)}(0)= \mathcal{C}^{(n)}(0)$. For odd $t$, the set
$\mathcal{D}^{(n)}(t)$ is then constructed by choosing,
independently for each individual in $\mathcal{D}^{(n)}(t-1)$, a
$\text{Binomial}(m,\gamma/n)$ distributed number of distinct
groups in $\mathcal{A}$, and, likewise, for even $t$, we construct
$\mathcal{D}^{(n)}(t)$ by choosing, independently for each group
in $\mathcal{D}^{(n)}(t-1)$, a $\text{Binomial}(n,\gamma/n)$
distributed number of distinct individuals in $\mathcal{V}$. Let
$X^{(n)}$ be a compound binomial random variable with generating
function
$$
g(s)=\E\left[s^{X^{(n)}}\right]=\left(1-\frac{\gamma}{n}+
\frac{\gamma}{n}\left(1-\frac{\gamma}{n}+\frac{\gamma}{n}s\right)^n\right)^m,
$$
and let $\{X^{(n)}(t): t\geq 0\}$ be a branching process with
offspring distribution $X^{(n)}$. Furthermore, write
$Y^{(n)}(t)=\sum_{s=0}^t X^{(n)}(s)$ for the total progeny of the
branching process at time $t$. Then, for even $t$, the number of
individuals (not necessarily distinct) that have been chosen in
the construction of the process $\mathcal{C}^{(n)}(t)$, has the
same distribution as $Y^{(n)}(t/2)$, and the number of groups (not
necessarily distinct) that have been chosen is (stochastically)
strictly smaller than $Y^{(n)}(t/2)$.

Define $\mu_n:=\mathbf{E}[X^{(n)}]=\gamma^2m/n$ and note that
$\beta\gamma^2(1-\frac1n) \leq \mu_n\leq \beta\gamma^2$, so
$\mu_n\to\mu:=\beta\gamma^2$ as $n\to\infty$. By a well known
result in branching process theory, we have that
$\mu_n^{-t}X^{(n)}(t)\to W^{(n)}$ almost surely as $t\to\infty$,
where $W^{(n)}$ is a random variable with $W^{(n)}\equiv 0$ if and
only if $\mu_n\leq 1$; see Athreya and Ney (1972). Furthermore,
$W^{(n)}\to W$ in distribution as $n\to\infty$, where $W$ is the
corresponding limiting random variable for the branching process
$\{X(t):t\geq 0\}$ with offspring generating function
$\mathbf{E}[s^{X(1)}]=\exp\{\beta\gamma(e^{\gamma(s-1)}-1)\}$.
Thus
$$
X^{(n)}(t)=\mu_n^t(W^{(n)}+o_t(1))=\mu_n^t(W+o_n(1)+o_t(1))\leq
\mu^t(O_n(1)+O_t(1)),
$$
where $o_x(\cdot)$ and $O_x(\cdot)$ is the usual order notation
when $x\to\infty$. It follows that $Y^{(n)}(t)\leq
\mu^t(O_n(1)+O_t(t))$, and, when we set $t=\lfloor \kappa \log n
\rfloor$ with $1/\kappa>\max\{2\log\mu,0\}$, we get
$$
Y^{(n)}(\lfloor \kappa \log n \rfloor) \leq e^{\log \mu \lfloor
\kappa \log n \rfloor}O_n(\log n)=o_n(\sqrt{n}).
$$
Now note that, if all individuals and groups that have been chosen
in the construction of $\mathcal{C}^{(n)}(t)$ are distinct, then
clearly the subgraph of $\mathcal{B}^{(n)}$ induced by
$\mathcal{C}^{(n)}(t)$ is a tree. Thus the probability in the
statement of the lemma is greater than
\begin{align*}
\prod_{k=1}^{Y^{(n)}(\lfloor\kappa\log n\rfloor)}\left(1-\frac{k}{n}\right) &
\left(1-\frac{k}{m}\right)=\\
&=\exp\bigg\{\sum_{k=1}^{o(\sqrt{n})}\left(\log\left(1-\frac{k}{n}\right)+
\log\left(1-\frac{k}{m}\right)\right)\bigg\}\\
&=\exp\bigg\{-\left(1+\frac1\beta\right)\sum_{k=1}^{o(\sqrt{n})}
\left(\frac{k}{n}+O\left(\frac{k^2}{n^2}\right)\right)\bigg\}\\
&=\exp\{o(1)\}\to 1,
\end{align*}
and the lemma is proved.\hfill$\Box$\medskip

\noindent\textbf{Proof of Theorem \ref{main}} The idea of the
proof is to construct a branching process $\{Z^{(n)}(t):t\geq 0\}$
that counts the number of individuals infected by the epidemic in
its initial stage, though not necessarily in chronological order.
The branching process will be defined in such a way that, if it
goes extinct, then its total progeny will be equal to the final
size of the epidemic, while, if it explodes, the epidemic will
have infected a number of individuals which is increasing
polynomially in $n$. As $n\to\infty$, we will have $Z^{(n)}\to Z$
--- where $Z$ is the branching process in the formulation of the
theorem --- and the theorem thus follows.

First, we describe the initial spread of the disease among the
individuals/groups in the set $\mathcal{C}^{(n)}(\lfloor
\kappa\log n\rfloor)$ with a process $\{\mathcal{E}^{(n)}(t):0\leq
t \leq \lfloor \kappa\log n\rfloor\}$. We only consider the case
when the subgraph of $\mathcal{B}^{(n)}$ induced by
$\mathcal{C}^{(n)}(\lfloor \kappa\log n\rfloor)$ is a tree. By
Lemma \ref{fodelsedag}, the probability of the complimentary set
tends to zero as $n\to\infty$, so we can disregard it.
Furthermore, our construction will be such that
$\mathcal{E}^{(n)}(t)\subseteq \mathcal{C}^{(n)}(t)$ for all $t$,
implying that nodes of $\mathcal{E}^{(n)}(\lfloor \kappa\log
n\rfloor)$ constitute a tree itself if seen as a subgraph of
$\mathcal{B}^{(n)}$.

The construction of $\mathcal{E}^{(n)}(t)$ is similar to the
construction of $\mathcal{C}^{(n)}(t)$ described in Lemma
\ref{fodelsedag}. Namely, we will define a sequence
$\{\mathcal{F}^{(n)}(t):0\leq t \leq \lfloor \kappa \log n
\rfloor\}$ and then set $\mathcal{E}^{(n)}(t)=\cup_{0\leq s \leq
t}\mathcal{F}^{(n)}(s)$. To this end, first let
$\mathcal{F}^{(n)}(0)$ consist of the initial infective, that is
$\mathcal{F}^{(n)}(0)=\{v_1\}$. Then, for odd $t$, let
$\mathcal{F}^{(n)}(t)$ consist of the groups of the individuals in
$\mathcal{F}^{(n)}(t-1)$ that are at distance $t$ from $v_1$, that
is,
$$
\mathcal{F}^{(n)}(t)=\{a\in\mathcal{D}^{(n)}(t):\exists v\in\mathcal{F}^{(n)}(t-1)
\textrm{ such that }v\in a\}.
$$
To define $\mathcal{F}^{(n)}(t)$ for even $t$, recall the
percolation representation of the set of ultimately infected
individuals in $\mathcal{G}^{(n)}$ described before the theorem.
For two vertices $v,w$ belonging to a group $a$, write
$v\overset{a}{\leftrightarrow}w$ for the event that there exists a
path of open edges --- that is, edges that can be used for disease
transmission --- connecting $v$ and $w$, with the additional
property that \emph{the whole path is contained in group $a$}.
Furthermore, let $\mathcal{K}_{a,v}=\{w\in
a:v\overset{a}{\leftrightarrow} w\}$. This is to be thought of as
the local outbreak in group $a$ caused by individual $v$, if $v$
itself becomes infected from outside of group $a$. As pointed out
before the formulation of the theorem, given that $|a|=k$, we have
that $|\mathcal{K}_{a,v}|\sim F_k$, where $F_k$ is the
distribution of the final size of of a homogeneous Reed-Frost
epidemic initiated by a single individual in a population of size
$k$. Note that $|a|\sim\text{Binomial}(n,\gamma/n)$, and, for
future use, let $R^{(n)}$ be a random variable with distribution
$\sum_kF_kP(|a|=k)$, that is, the size of a local outbreak in a
group, not conditioning on the group size. Now, for even $t$,
define $\mathcal{F}^{(n)}(t)$ to be the individuals infected in
the local outbreaks caused by the individuals in
$\mathcal{F}^{(n)}(t-2)$, that is,
$$
\mathcal{F}^{(n)}(t)=\{w\in \mathcal{K}_{a,v}: a\in\mathcal{F}^{(n)}(t-1),
v\in\mathcal{F}^{(n)}(t-2)\}.
$$
We will now study the growth of $|\mathcal{E}^{(n)}(t)|$. To this
end, for $0\leq t \leq \frac12\lfloor \kappa \log n \rfloor$,
define $Z^{(n)}(t)=|\mathcal{F}^{(n)}(2t)|$. Then, since the
subgraph of $\mathcal{B}^{(n)}$ induced by $\mathcal{C}^{(n)}(2t)$
is a tree for $t\leq \frac12\lfloor \kappa \log n \rfloor$, by
construction, $Z^{(n)}(t)$ is a branching process with a compound
binomial offspring distribution. The generating function of the
offspring distribution is
\begin{equation}\label{eq:f_n}
f_n(s)=\E\left[s^{Z^{(n)}(1)}\right]=\left(1-\frac{\gamma}{n}+
\frac{\gamma}{n}\E\left[s^{R^{(n)}}\right]\right)^m.
\end{equation}
For $t\geq\frac12\lfloor \kappa \log n \rfloor$ we let
$Z^{(n)}(t)$ evolve by the same branching mechanism, that is, as a
discrete time branching process with offspring distribution
defined by (\ref{eq:f_n}). For $t\geq\frac12\lfloor \kappa \log n
\rfloor$ however, $Z^{(n)}(t)$ is no longer related to the
epidemic process.

Let $\rho_n$ be the smallest non-negative root of
$\rho_n=f_n(\rho_n)$, and define
$A_n=\{\lim_{t\to\infty}Z^{(n)}(t)=0\}$ (the extinction set) and
$T_n=\inf\{t:Z^{(n)}(t)=0\}$ (the extinction time). Then
$\P(A_n)=\rho_n$, and, since $T_n$ is finite on $A_n$, the total
progeny of the branching process, $\sum_{t=0}^{T_n}Z^{(n)}(t)$, is
finite on $A_n$. As $n\to\infty$, we have that $Z^{(n)}\to Z$ in
distribution and $\rho_n\to \rho$, where $Z(t)$ is the branching
process appearing in the theorem and $\rho$ the equivalent of
$\rho_n$ in this process. Furthermore, $T^{(n)}\to T$ in
distribution, where $T=\inf\{t:Z(t)=0\}$. It follows that
$E=\sum_{t=0}^\infty Z(t)$ is finite on the event
$A=\{\lim_{t\to\infty}Z(t)=0\}$. This implies that, on $A_n$, the
final size  $E^{(n)}$ of the epidemic converges to $E$ as
$n\to\infty$. On the complementary sets $A_n^c$, the process
$Z^{(n)}(t)$ grows exponentially in $t$. More precisely,
$$
Z^{(n)}(\lfloor\kappa\log n\rfloor)=n^{\kappa\log R_0}(W'+o(1)),
$$
where $W'>0$ due to the conditioning on explosion (recall the
proof of Lemma \ref{fodelsedag}). Furthermore, we have that
$Z^{(n)}(\lfloor\kappa\log n\rfloor)\leq E^{(n)}$ on $A_n^c$, and
thus $E^{(n)}\to\infty$ on $A_n^c$. This proves the theorem.
\hfill$\Box$

\section{The final outcome of the epidemic}

The branching process approach from the previous section gives
basically no information on the behavior of the epidemic in the
case of explosion. In this section we will elaborate a bit on this
problem.

As already described, one way of getting a grip of the final
outcome of the epidemic, is to consider an edge percolation
process on the underlying graph, where each edge in the graph is
independently removed with probability $1-p$ and kept with
probability $p$. The vertices that belong to the component of the
initial infective in the graph so obtained correspond to the
individuals that have experienced the infection at the end of the
epidemic. If the structure of the thinned graph is known, then
this observation might be useful in investigating the final size
of the epidemic. For instance, if there is a unique giant
component in the thinned graph --- that is, if the outcome of the
percolation process contains a unique cluster of order $n$ ---
then the relative size of this component gives the probability of
an outbreak of order $n$ in the epidemic. Such an outbreak is
often referred to as a major outbreak, and, in most epidemic
models, the probability of such an outbreak coincides with the
probability of explosion in the branching process describing the
initial stages of the epidemic (denoted by $\pi$ in this paper).
This however require additional arguments.

In our case, the social network is a random intersection graph
with $\alpha=1$. Unfortunately, to date there are no rigorous
results concerning the component structure in a random
intersection graph with $\alpha=1$, but see Behrisch (2007) for
results when $\alpha\neq 1$. Also, in Newman (2003:2), (implicit)
expressions for the size of the largest component in a random
graph construction which is similar to the random intersection
graph are derived by heuristic means and it is observed that the
relative final size of the giant component seems to decrease as
the clustering in the graph increases. An argument in support of
the claim that high clustering in a graph causes the components to
be small is the following: Consider an arbitrary graph with $n$
vertices and $k=O(n)$ edges and assume that the clustering equals
1. This implies that all subgraphs are complete. Hence, with
$n_{max}$ denoting the size of the largest subgraph, we have that
the number of edges in the maximal subgraph is $n_{max}\choose 2$.
It follows that $n_{max}\leq O(\sqrt{k})=O(\sqrt{n})$, that is,
the relative size of the largest component tends to zero,

Indeed, the lack of rigorous results concerning the components in
a random intersection graph with $\alpha=1$ makes it harder to
study the final size of an epidemic on such a graph. A second
complicating circumstance is that thinning a random intersection
graph gives rise to a graph that no longer belongs to the class of
random intersection graphs; see the below proposition. This means
that, even if there would be results for the component structure,
these would not be applicable to a thinned graph. Hence it remains
an open problem to quantify the final outcome of the epidemic.

\begin{prop}
Let $\Theta_{p} (\mathcal{G}_{\beta, \gamma}^{(n)})$ denote the graph generated by removing edges in $\mathcal{G}_{\beta, \gamma}^{(n)}$
independently with probability $1-p$. It does not exist $\beta^{\prime} = \beta^{\prime}(\beta, \gamma, p)$
and $\gamma^{\prime} = \gamma^{\prime}(\beta, \gamma, p)$ such
that $\Theta_{p} (\mathcal{G}_{\beta, \gamma}^{(n)})
\stackrel{d}{=} \mathcal{G}_{\beta^{\prime},
\gamma^{\prime}}^{(n)}$ for every $n$.
\end{prop}

\noindent \textbf{Proof.} The idea of the proof is to observe that
certain types of subgraphs will appear with different frequency in
$\Theta_{p} (\mathcal{G}_{\beta, \gamma}^{(n)})$ as compared to
$\mathcal{G}_{\beta, \gamma}^{(n)}$. The subgraph that we will
consider consists of four vertices and five edges:
$$
\xymatrix{
{\bullet} \ar@{-}[d] \ar@{-}[dr] \ar@{-}[r] & {\bullet} \ar@{-}[d]\\
{\bullet} \ar@{-}[r] & {\bullet} }
$$
Write $K_4'$ for this graph type, and note that it can be obtained
for instance by removing one edge from a complete subgraph with
four vertices, a graph type that we denote by $K_4$. Furthermore,
we introduce the term \emph{vertex-induced} subgraph, for a
subgraph of some given graph such that the subgraph consists of a
subset of the vertices in the original graph together with all
edges between these vertices that are present in the original
graph.

The number $X_4$ of vertex-induced subgraphs of type $K_{4}$ in
the random intersection graph $\mathcal{G}_{\beta, \gamma}^{(n)}$
dominates the number of groups of size four in the construction of
the random intersection graph. Since the size of a fixed group is
${\rm Binomial}(n,\gamma/n)$ distributed, the number of groups of size four is ${\rm Binomial}(\lfloor \beta n\rfloor,
\binom{n}{4}(\gamma/n)^{4}(1-\gamma/n)^{n-4})$ distributed, and hence $\E[X_4]
\geq O(n)$. It follows that the number $X_4'(p)$ of vertex-induced
subgraphs of type $K_4'$ in the thinned graph $\Theta_{p}
(\mathcal{G}_{\beta, \gamma}^{(n)})$ is also at least of the order
$n$, since, as mentioned, one way of obtaining graphs of type
$K_4'$ is to remove one edge in graphs of type $K_{4}$, that is,
$\E[X_4'(p)] \ge (1-p)p^{5}\E[X_4] \geq O(n)$.

Now consider the number $X_4'$ of vertex-induced subgraphs of type
$K_4'$ in the random intersection graph $\mathcal{G}_{\beta,
\gamma}^{(n)}$. This number is related to the number of ways that
four individuals $v_1,\ldots,v_4$ can be assigned to different
groups so that a graph of type $K_4'$ is obtained. Consider for
instance the following graph of type $K_4'$:
$$
\xymatrix{ {}^{v_{1}}\hspace{-18mm} & {\bullet} \ar@{-}[d]
\ar@{-}[dr] \ar@{-}[r] &
{\bullet} \ar@{-}[d] & \hspace{-18mm} {}^{v_{2}}\\
{}_{v_{3}}\hspace{-18mm} & {\bullet} \ar@{-}[r] & {\bullet} &
\hspace{-18mm} {}_{v_{4}} }
$$
Write $\{(v_{i_1},v_{i_2})(v_{j_1},v_{j_2})(v_{k_1},v_{k_2})\}$
for the event that $v_{i_1}$ and $v_{i_2}$ share a group, that
$v_{j_1}$ and $v_{j_2}$ share \emph{another} group and that
$v_{k_1}$ and $v_{k_2}$ share yet another group. Then, for the
above graph to arise, the individuals $v_2$ and $v_3$ cannot share any group --- the probability that they avoid doing so goes to 1 as $n\to\infty$ --- and, in addition, one of the following events must occur:
$$
\begin{array}{ll}
\{(v_{1}, v_{2}, v_{4})(v_{1}, v_{3}, v_{4})\}\\
\{(v_{1}, v_{2}, v_{4})(v_{1}, v_{3})(v_{3},
v_{4})\}\\
\{(v_{1}, v_{3}, v_{4})(v_{1}, v_{2})(v_{2},
v_{4})\}\\
\{(v_{1}, v_{2})(v_{1}, v_{4})(v_{2}, v_{4})(v_{1}, v_{3}) (v_{3},
v_{4})\}.
\end{array}
$$
It follows that
\begin{align*}
\E[X_4'] &\le n^4\binom{\lfloor \beta n \rfloor}{2}(\gamma/n)^{6}
+ 2n^4\binom{\lfloor
\beta n \rfloor}{3}(\gamma/n)^{7} \\
&+ n^4\binom{\lfloor \beta n \rfloor}{5}(\gamma/n)^{10} = O(1).
\end{align*}

The number of vertex-induced subgraphs of type $K_4'$ in a random
intersection graph is hence finite, while, in a thinned random
intersection graph, it is of order $n$. This proves the
proposition. \hfill$\Box$

\section{Numerical results}

\begin{figure}
    \begin{center}
        \begin{minipage}[c]{0.90\textwidth}
            \begin{center}
                \includegraphics[width=1\textwidth, height=0.43\textheight]{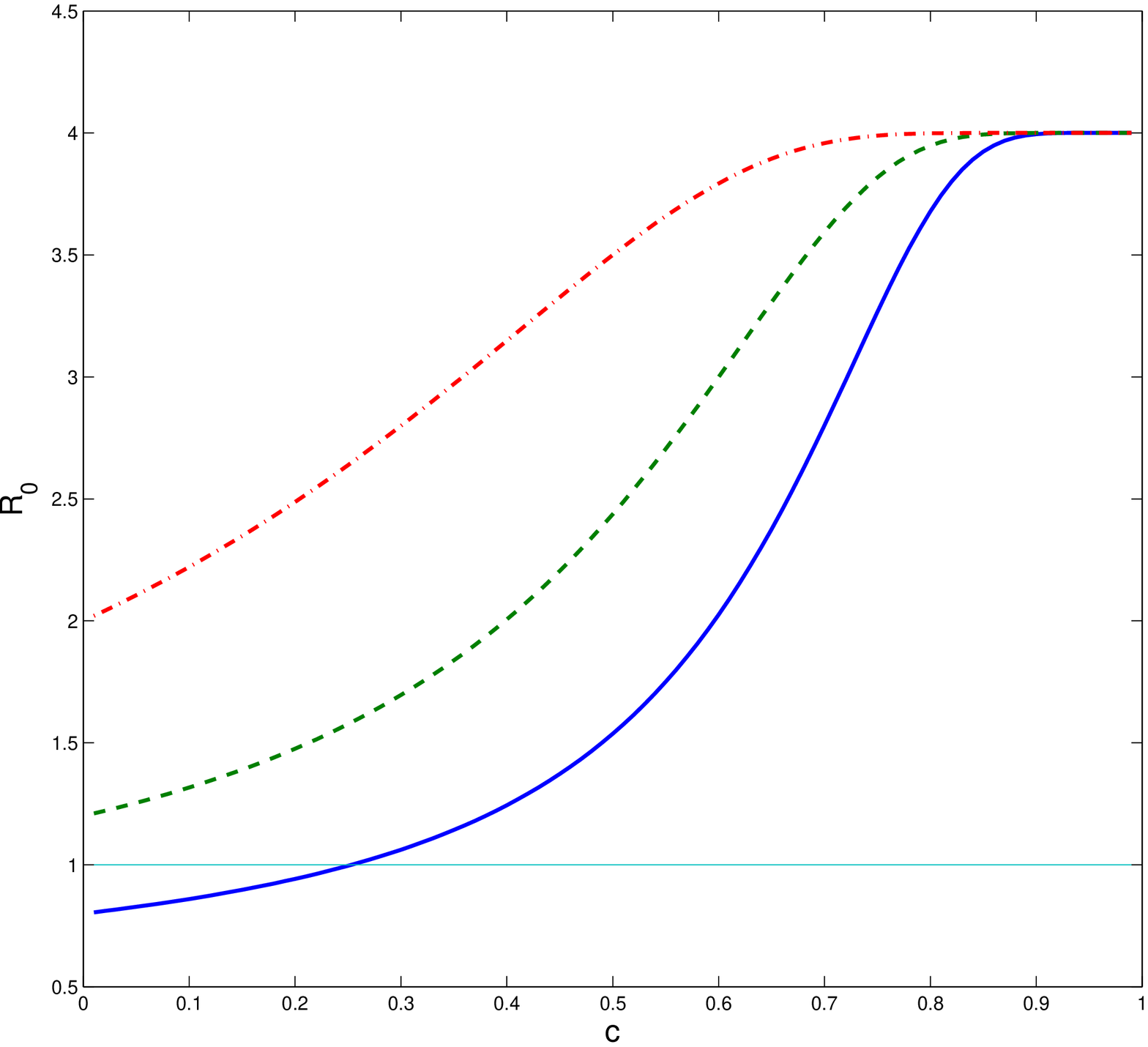} \qquad
            \end{center}
        \end{minipage} \qquad
        \begin{minipage}[c]{0.90\textwidth}
            \begin{center}
                \includegraphics[width=1\textwidth, height=0.43\textheight]{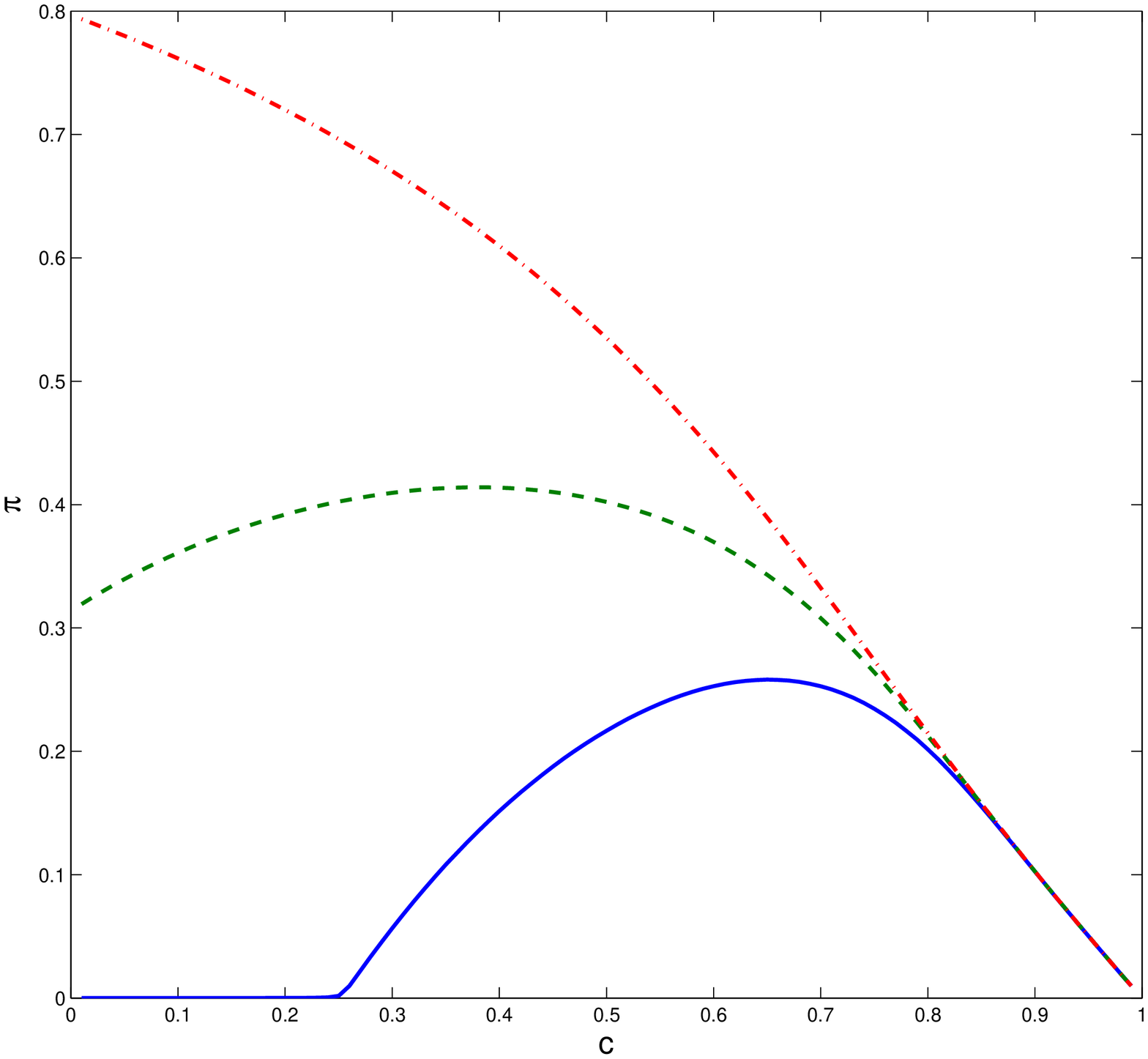}
            \end{center}
        \end{minipage}
    \end{center}
    \caption{In the top figure, $R_0$ is plotted against $c$ for fixed $\mu=4$ and
    for three choices of $p$\,: $p = 0.2$
    (------), $0.3$ (-- -- --), and $0.5$ (-- $\cdot$ --). The bottom
    figure shows how the probability $\pi$ of explosion in the
    epidemic varies with $c$, for fixed $\mu = 4$, and the same choices
    of $p$ as above.}\label{figuren}
\end{figure}

In this section, we investigate numerically the epidemic threshold
$R_0$ and the probability $\pi$ of explosion in the epidemic ---
recall Corollary \ref{large_outbreak}. We have $R_0=\beta
\gamma\E[R]$, where the distribution of $R$ is specified in
(\ref{eq:R}), and $\pi:=1-\rho$, where $\rho$ is the smallest
non-negative root of the equation $f(\rho)=\rho$ and $f$ is the
generating function specified in (\ref{eq:f}). Using the
recursive formulas for the distribution $F_k$ of the final size a
Reed-Frost epidemic in a homogeneous population of size $k$ ---
see e.g.\ Andersson and Britton (2000: Section 1.2)
--- numerical values of $R_0$ and $\pi$ are easily obtained for
fixed values of $\beta$ and $\gamma$.

We are particularly interested in how $\rho$ and $R_{0}$ are
affected when the (asymptotic) clustering $c=(1+\beta\gamma)^{-1}$
is varied, and, to be able to compare results for different values
of $c$, the mean degree $\mu=\beta\gamma^2$ in the graph
is kept fixed. In Figure 1, the parameters $R_0$ and $\pi$ are
plotted against $c$ for three different values of the infection
probability $p$. Let us comment a bit on these plots.

First we investigate the value of $R_0$ in the limit as $c\to 0$.
Since $\mu$ is fixed, we have that $c=(1+\mu/\gamma)^{-1}\to 0$
implies that $\gamma\to 0$ as well. Asymptotically, the degree
distribution in our random intersection graph is compound Poisson
with generating function
$$
g(s) = e^{\beta\gamma \left(e^{\gamma(s-1)} -
1\right)}
$$
which converges to $e^{\mu(s-1)}$ as $\gamma\to
0$. The limiting degree distribution as $c\to 0$ is hence
Poisson($\mu$), and, after thinning the graph, removing each edge
independently with probability $1-p$, the degrees are
Poisson distributed with mean $p\mu$. Since the graph obtained by
such a thinning can be thought of as representing the outcome of
the epidemic, it is reasonable to suspect that $R_0=p\mu$ in the
limit as $c\to 0$. Indeed, it can be seen in the top plot in
Figure \ref{figuren} that $R_{0} \to p\mu$ as $c \to 0$.

In the top plot in Figure 1, it can also be seen that $R_0$
increases with $c$, that is, higher clustering makes it easier for
epidemics to take off. This is in line with findings in Newman
(2003:2). Let us give a heuristic explanation of why this should be the
case: First note that, since the mean degree $\mu$ in the graph is
fixed, an increase in $c=(1+\mu/\gamma)^{-1}$ is equivalent to an
increase in $\gamma$ and a decrease in $\beta$ of the order
$\gamma^{-2}$. Also, recall that the mean number of groups that an
individual is a member of is $\beta\gamma$ and the mean group size
is $\gamma$. Hence, increased clustering with fixed mean degree
means that individuals are members of fewer but larger groups.
Combining this with the observation that the probability for an
individual to avoid infection from some index case with whom
he/she shares a group decreases geometrically with the group size,
it follows that it should be easier for an epidemic to take off
when the clustering is large. In fact, we have that $R_{0}\to\mu$
as $c \to 1$, that is, in the limit of large clustering, the
infection probability $p$ does not matter (as long as it is
positive) for the value of $R_0$.

The bottom plot in Figure 1 shows how the probability $\pi$ of
explosion in the epidemic varies with $c$. For instance, it can be
seen that $\pi\to 0$ as $c\to 1$. In Section 4, we argued that the
relative size of the largest component in a graph with maximal
clustering is 0 in the limit of large graph size. If the
probability of explosion in the epidemic coincides with the
relative size of the largest component in the graph representing
the outcome of the epidemic, then indeed it follows from this that
$\pi\to 0$ as $c\to 1$. Furthermore, it is interesting to note
that the decrease of $\pi$ towards 0 is not monotone for all
values of $p$. Clearly, if a low value of $c$ prevents explosion,
while explosion is possible for a larger value of $c$
--- this is the case for instance for $p=0.2$ --- then we will see an
increase in $\pi$ from 0 to a positive value when the threshold is
passed. But, as the curve for $p=0.3$ reveals, even if $\pi$ is
positive already at $c=0$, it can be the case that it increases
with $c$ in some interval before it starts to decrease.

\section{Discussion}

In the present paper, we have analyzed how the clustering in a
random network affects how an infectious disease propagates in the
network, assuming the size of the network to be large.  In
particular, using a random intersection graph construction, we
have rigorously derived the limiting probability of an explosion
in the epidemic and a threshold parameter indicating if this
probability is 0 or positive.

The motivation for analyzing an epidemic on a random network with
positive clustering is of course that most empirical social
networks manifest positive clustering, so predictions based on
epidemic models neglecting such clustering, i.e.\ most epidemic
models, must be interpreted with caution. There are of course
several other features in empirical networks, not considered in
the present paper, that should also be taken into account for
predictions to be reliable. One such feature is the degree
distribution, which in many social networks has been observed to follow a power-law distribution. The graph model used in this
paper gives compound Poisson distributions for the degrees, but
the model is generalized in Deijfen and Kets (2007) to allow for
power-law degree distributions. It would be interesting to study
how an epidemic on such a generalized graph is affected by the
exponent in the power-law. Another feature that has been observed
in many social networks is positive degree correlation, that is,
individuals with high (low) degree tend to be connected to other
individuals with high (low) degree. Because of the group
structure, this is likely to be the case in a random intersection
graph, but it remains to quantify the correlation.

A possible generalization of the studied model would be to
distinguish between different types of individuals, and to assume
that both network characteristics as well as transmission
probabilities depend on the type of an individual; see e.g.\ Ball
and Clancy (1993). Another extension, motivated by real-world
epidemics, is to leave the Reed-Frost paradigm, in which the
events that different neighbors of a given infective becomes
infected are independent. If for example the infectious period is
taken to be random, then these events are positively correlated;
see e.g.\ Andersson and Britton (2000). Unfortunately, by relaxing
the independence assumption the analysis of the model becomes much more
complicated.

Perhaps the most obvious continuation of the present work is
however to derive fully rigorous results about the final size of
the epidemic in case of explosion. The (relative) final
\emph{size} of the epidemic then most likely coincides with the
\emph{probability} of explosion, a quantity derived in the present
paper, but a proof of this is still missing.

\section*{References}

\noindent Andersson, H. (1998): Limit theorems for a random graph
epidemic model, \emph{Annals of Applied Probability} \textbf{8},
1331--1349.\medskip

\noindent Andersson, H. (1999): Epidemic models and social networks,
\emph{The Mathematical Scientist} \textbf{24}, 128--147.\medskip

\noindent Andersson, H. and Britton, T. (2000): \emph{Epidemic
models and their statistical analysis}, Springer.\medskip

\noindent Athreya, K.\ B.\ and Ney, P.\ E.\ (1972):
\emph{Branching processes}, Springer.\medskip

\noindent von Bahr, B. and Martin-L\"{o}f, A. (1980): Threshold
limit theorems for some epidemic processes, \emph{Advances in
Applied Probability} \textbf{12}, 319--349.\medskip

\noindent Ball, F. and Clancy, D. (1993): The final size and
severity of a generalized stochastic multi-type epidemic model,
\emph{Advances in Applied Probability} \textbf{25},
721--736.\medskip

\noindent Ball, F., Mollison, D. and Scalia-Tomba, G. (1997):
Epidemics with two levels of mixing, \emph{Annals of Applied
Probability} \textbf{7}, 46--89.\medskip

\noindent Behrisch, M. (2007): Component evolution in random
intersection graphs, \emph{Electronic Journal of Combinatorics}
\textbf{14}.\medskip

\noindent Deijfen, M. and Kets, W. (2007): Random intersection
graphs with tunable degree distribution and clustering, preprint,
available at www.math.su.se/$\sim$mia.\medskip

\noindent Dorogovtsev, S. and Mendes, J. (2003): \emph{Evolution of
Networks, from Biological Nets to the Internet and WWW}, Oxford
University Press.\medskip

\noindent Durrett, R. (2006): \emph{Random graph dynamics},
Cambridge University Press.\medskip

\noindent Fill, J., Scheinerman, E. and Singer-Cohen, K. (2000):
Random intersection graphs when $m=\omega(n)$: an equivalence
theorem relating the evolution of the $G(n,m,p)$ and $G(n,p)$
models, \emph{Random Structures $\&$ Algorithms} \textbf{16},
156--176.\medskip

\noindent Godehardt, E. and Jaworski, J. (2002): Two models of
random intersection graphs for classification, in \emph{Exploratory
data analysis in empirical research}, eds Schwaiger M. and Opitz,
O., Springer, 67--81.\medskip

\noindent Karo\'{n}ski, M., Scheinerman, E. and Singer-Cohen, K.
(1999): On random intersection graphs: the subgraphs problem,
\emph{Combinatorics Probability $\&$ Computing} \textbf{8},
131--159.\medskip

\noindent Neal, P. (2004): SIR epidemics on Bernoulli graphs,
\emph{Journal of Applied Probability}, \textbf{40}, 779--782.\medskip

\noindent Neal, P. (2006): Multitype randomized Reed-Frost
epidemics and epidemics upon graphs, \emph{Annals of Applied
Probability} \textbf{16}, 1166--1189.\medskip

\noindent Newman, M. (2003:1): The structure and function of
complex networks, \emph{SIAM Reviews} \textbf{45},
167--256.\medskip

\noindent Newman, M. (2003:2): Properties of highly clustered
networks, \emph{Physical Review E} \textbf{68}, 026121.\medskip

\noindent Newman, M., Barab\'{a}si, A. and Watts, D. (2006):
\emph{The structure and dynamics of networks}, Princeton University
Press.\medskip

\noindent Singer, K. (1995): \emph{Random intersection graphs}, PhD
thesis, Johns Hopkins University.\medskip

\noindent Stark, D. (2004): The vertex degree distribution of random
intersection graphs, \emph{Random Structures $\&$ Algorithms}
\textbf{24}, 249--258.\medskip

\noindent Trapman, P. (2007): On analytical approaches to epidemics
on networks, \emph{Theoretical Population Biology} \textbf{71},
160--173.

\end{document}